\documentstyle[amscd,amssymb,amsthm,verbatim,amsopn,amsopn,12pt]{amsart}

\theoremstyle{plain}
\newtheorem{Thm}{Theorem}[section]

\newtheorem{Prop}[Thm]{Proposition}
\newtheorem{Cor}[Thm]{Corollary}
\theoremstyle{definition}

\newtheorem{Conj}[Thm]{Conjecture}

\newcommand{\bnum}{\begin{enumerate}}
\newcommand{\enum}{\end{enumerate}}

\DeclareMathOperator{\girth}{girth}
\DeclareMathOperator{\diam}{diam}
\newcommand{\Z}{\mathbb{Z}} 
\newcommand{\N}{\mathbb{N}}

\begin{document}
\title[Relative non-commuting graph of a finite ring]%
	{Relative non-commuting graph of a finite ring}
\author[J. Dutta and D. K. Basnet]%
	{Jutirekha Dutta \&  Dhiren K. Basnet*}
\thanks{*Corresponding author}
\date{}
\maketitle
\begin{center}\small{\it Department of Mathematical Sciences, Tezpur University,}
\end{center}
\begin{center}\small{\it Napaam-784028, Sonitpur, Assam, India.}
\end{center}

\begin{center}\small{\it Email: jutirekhadutta@@yahoo.com  and dbasnet@@tezu.ernet.in}
\end{center}

\medskip

\noindent \textit{\small{\textbf{Abstract:}
 Let $S$ be a subring of a finite ring $R$ and $C_R(S) = \{r \in R : rs = sr \;\forall\; s \in S\}$. The relative non-commuting graph of the subring $S$ in $R$, denoted by $\Gamma_{S, R}$, is a simple undirected graph whose vertex set is $R \setminus C_R(S)$ and two distinct vertices $a, b$ are adjacent if and only if $a$ or $b \in S$ and $ab \neq ba$. In this paper, we discuss some properties of $\Gamma_{S, R}$, determine  diameter, girth, some dominating sets and chromatic index for $\Gamma_{S, R}$. Also, we derive some connections between  $\Gamma_{S, R}$ and the relative commuting probability of $S$ in $R$. Finally, we show that  the relative non-commuting graphs of two relative $\Z$-isoclinic pairs of rings are isomorphic under some conditions.}}

\bigskip

 

\medskip

\section{Introduction}
Let $R$ be a finite ring with subring  $S$. Let $C_R(S) = \{r \in R : rs = sr \;\forall\; s \in S\}$. The relative non-commuting graph of the subring $S$ in $R$, denoted by $\Gamma_{S, R}$, is defined as a simple undirected graph whose vertex set is $R \setminus C_R(S)$ and two distinct vertices $a, b$ are adjacent if and only if $a$ or $b \in S$ and $ab \neq ba$. For $S = R$, we have $\Gamma_{S, R} = \Gamma_R$, the non-commuting graph of $R$.  The notion of non-commuting graph of a finite ring  was introduced by Erfanian et al. \cite{ekn15} in the year 2015. The study of algebraic structures by means of graph theoretical properties became more popular during the last decade
(see \cite{abdollahi07, abdollahi08, abdollahiakbarimaimani06, beck88, ov11} etc.). Motivated by the works of Erfanian et al. \cite{TE13}, in this paper, we  obtain some graphs that are not isomorphic to   
$\Gamma_{S, R}$ for any ring $R$ with  subring $S$. We also 
determine  diameter, girth, some dominating sets and chromatic index for $\Gamma_{S, R}$ and derive some connections between  $\Gamma_{S, R}$ and the relative commuting probability of $S$ in $R$.
Recall that the relative commuting probability of a subring $S$ in a finite ring $R$, denoted by $\Pr(S, R)$,  is the probability that a randomly chosen pair of elements, one from $S$ and the other from $R$ commute.  That is
\[
\Pr(S, R) = \frac{|\{(s, r) \in S \times R : sr = rs\}|}{|S||R|}.
\]
This notion was introduced and studied in \cite{duttabasnetnath}. Note that  $\Pr(R, R)$ is the commuting probability of $R$, a notion introduced by MacHale \cite{machale}. In the last section, we show that  the relative non-commuting graphs of two relative $\Z$-isoclinic pairs of rings are isomorphic under some conditions.

For a graph $G$, we write $V(G)$ and $E(G)$ to denote the set of vertices and the set of edges of $G$ respectively. We write $\deg(v)$ to denote  the degree of a vertex $v$, which is the number of edges incident on $v$. Let $\diam(G)$ and $\girth(G)$  be the diameter and girth of a graph $G$ respectively. Recall that 
$\diam(G) = \max\{d(x, y) : x, y \in V(G)\}$, where $d(x, y)$ is the length of the shortest path from $x$ to $y$; and $\girth(G)$ is the length of the shortest cycle obtained in $G$.  
A graph $G$ is called connected if there is a path between every pair of vertices. A star graph is a tree on $n$ vertices in which one vertex has degree $n - 1$ and the others have degree $1$. A bipartite graph is a graph whose vertex set can be partitioned into two disjoint parts in such a way that the two end vertices of every edge lie in different parts. A complete bipartite graph is a bipartite graph such that two vertices are adjacent if and only if they lie in different parts. A complete graph is a graph in which every pair of distinct vertices is adjacent. Throughout the paper $R$ denotes a finite non-commutative ring.


\section{Some properties of $\Gamma_{S, R}$}
Let $S$ be a subring of a ring $R$,  $r \in R$ and $A \subseteq R$. We write $C_R(r) := \{x \in R : xr =rx\}$,  $C_S(r) := C_R(r) \cap S$ and $C_R(A) = \{x \in R : xa = ax \, \forall\, a \in A\}$. Note that $C_R(r)$ and $C_S(r)$ are subrings of $R$. Also 
$\underset{r \in R}{\cap} C_R(r) := Z(R)$ is  the center of $R$. We begin this section with the following useful result.

\begin{Prop}\label{deg-prop}
Let $S$ be a non-commutative subring of a ring $R$. Then
\begin{enumerate}
\item $\deg(r) = |R| - |C_R(r)|$ if $r \in V(\Gamma_{S, R}) \cap S$.
\item $\deg(r) = |S| - |C_S(r)|$ if $r \in V(\Gamma_{S, R}) \cap (R \setminus S)$.
\item $\Gamma_{S, R}$ is connected.
\item $\Gamma_{S, R}$ is empty graph if and only if $S$ is commutative.
\end{enumerate} 
\end{Prop}
\begin{pf}
The proof of part (a), (b) and (d) follow from the definition of $\Gamma_{S, R}$. For part (c), suppose $\Gamma_{S, R}$ has an isolated vertex, namely $v$. Then $\deg(v) = |R| - |C_R(v)| = 0$ or $|S| - |C_S(v)| = 0$ for $v \in S$ or $v \in R \setminus S$. Thus, in both cases $v \in C_R(S)$, a contradiction.
\end{pf}

In the following theorems we shall show that if $G$ is a star graph or an $n$-regular graph, where $n$ is a square free odd positive integer, then $G$ can not be realized by $\Gamma_{S, R}$ for any subring $S$ of a ring $R$. Also, $\Gamma_{S, R}$  is not a bipartite graph, for any proper subring $S$ of a ring $R$.

\begin{Thm}\label{star}
Let $S$ be a non-commutative subring of a ring $R$. Then $\Gamma_{S, R}$ is not a star graph.
\end{Thm}
\begin{pf}
Suppose, $\Gamma_{S, R}$ is a star graph, where $S$ is a non-commutative subring of $R$. 
Then all but one vertices of $\Gamma_{S, R}$ have degree $1$. 
Let $v$ be a vertex of  $\Gamma_{S, R}$ having degree $1$. Then, by Proposition \ref{deg-prop}, we have  $[R : C_R(v)] = |R|/(|R| - 1)$ or $[S : C_S(v)] = |S|/(|S| - 1)$ according as 
$v \in S$ or $v \in R \setminus S$; which is absurd.
Hence the result follows.
\end{pf}

\begin{Thm}\label{complete relative non-commuting}
Let $S$ be a proper non-commutative subring of a ring $R$. Then $\Gamma_{S, R}$ is not bipartite.
\end{Thm}

\begin{pf}
Let $\Gamma_{S, R}$ be a bipartite graph.
Then, there exist two disjoint subsets $S_1$ and $S_2$ of $V(\Gamma_{S, R})$ such that $|S_1| + |S_2| = |R| - |C_R(S)|$. Therefore, $S \cap S_1 = \phi$ or $S \cap S_2 = \phi$. So, $S \subseteq S_2$ or $S \subseteq S_1$. Without loss of generality we may assume that  $S \subseteq S_1$. Then, for $v \in S_1$ we have $vs = sv$ for all $s \in S \setminus C_R(S)$. Thus, $v \in Z(S) \subseteq C_R(S)$, a contradiction. Hence, the theorem follows.
\end{pf}

\begin{Thm}\label{n-regular}
Let $S$ be a non-commutative subring of a ring $R$. Then $\Gamma_{S, R}$ is not an $n$-regular graph for  any square free odd positive integer $n$.
\end{Thm}
\begin{pf}
Let $\Gamma_{S, R}$ be an $n$-regular graph.
Suppose, $n = p_1p_2\dots p_m$, where $p_i$'s are distinct odd primes. If $v \in V(\Gamma_{S, R}) \cap S$  then, by Proposition \ref{deg-prop}, we have
\[
n = \deg(v) = |R| - |C_R(v)| = |C_R(v)|([R : C_R(v)] - 1).
\]
 Here, $|C_R(v)| \neq 1$, as $0, v \in C_R(v)$. Thus $|C_R(v)| = \underset{p_i \in Q}{\prod}p_i$ and $[R : C_R(r)] - 1 = \underset{p_j \in P\setminus Q}{\prod}p_j$, where $Q \subseteq \{p_1, p_2, \dots, p_m\} = P$. So, $|R| = \underset{p_i \in Q}{\prod}p_i(\underset{p_j \in P \setminus Q}{\prod}p_j + 1)$. If $r \in R \setminus S$ then, using similar argument, we have $|S| = \underset{p_i \in T}{\prod}p_i(\underset{p_j \in P \setminus T}{\prod}p_j + 1)$, where $T \subseteq  P$. So,   $\underset{p_i \in T - (T \cap Q)}{\prod}p_i(\underset{p_j \in P \setminus T}{\prod}p_j + 1)$ divides $\underset{p_j \in P \setminus Q}{\prod}p_j + 1$, which is not possible. Hence, the theorem follows.
\end{pf}



We conclude this section showing that a complete graph can not be realized by $\Gamma_{S, R}$ for a subring $S$ of a ring $R$ with unity. 
\begin{Thm}
Let $R$ be a ring with unity and $S$ a subring of $R$. Then $\Gamma_{S, R}$ is not complete. 
\end{Thm}
\begin{pf}
 Suppose that there exists a subring $S$ of $R$ with unity such that  $\Gamma_{S, R}$ is complete. Then, for any $s \in V(\Gamma_{S, R}) \cap S$ we have 
\[
\deg(s) = |V(\Gamma_{S, R})| - 1 = |R| - |C_R(S)| - 1.
\]
By Proposition \ref{deg-prop}, we have $|R| - |C_R(s)| = |R| - |C_R(S)| - 1$. This gives $|C_R(S)| = 1$ and $|C_R(s)| = 2$, which is not possible, since $R$ is a ring with  unity. Hence, the result follows.
\end{pf}

\section{Diameter, girth, dominating set and chromatic index}

In this section, we obtain diameter, girth, some dominating sets and chromatic index of the graph $\Gamma_{S, R}$. 
\begin{Thm}\label{diam and girth}
Let $S$ be a non-commutative subring of a ring $R$. If $Z(S) = \{0\}$ then $\diam(\Gamma_{S, R}) = 2$ and $\girth(\Gamma_{S, R}) = 3$. 
\end{Thm}
\begin{pf}
Suppose, $v_1$ and $v_2$ are two vertices of $\Gamma_{S, R}$  such that they are not adjacent. So, there exist vertices $s_1, s_2 \in S$ such that $v_1s_1 \neq s_1v_1$ and $v_2s_2 \neq s_2v_2$. If $v_2$ is adjacent to $s_1$ or $v_1$ is adjacent to $s_2$, then $d(v_1, v_2) = 2$. Suppose that both are not adjacent, that is $v_1s_2 = s_2v_1$ and $v_2s_1 = s_1v_2$. Then $s_1 + s_2$ is adjacent to $v_1$ and $v_2$, which give $d(v_1, v_2) = 2$. Therefore, $\diam(\Gamma_{S, R}) = 2$.

In order to determine $\girth(\Gamma_{S, R})$, suppose that  $v, s \in V(\Gamma_{S, R})$ where $s \in S$ and $v, s$ are adjacent. So, there exist $v_1, v_2 \in V(\Gamma_{S, R})$ such that $v$ and $s$ are adjacent to $v_1$ and $v_2$ respectively. If $v, v_2$ or $s, v_1$ are adjacent then $\{v, s, v_2\}$ or  $\{v, s, v_1\}$ is a cycle of length $3$ in $\Gamma_{S, R}$. If both are not adjacent then $v_1 + v_2$ is adjacent to $v$ and $s$. Therefore, $\{v, s, v_1 + v_2\}$ is a cycle of length $3$ in $\Gamma_{S, R}$. Hence, $\girth(\Gamma_{S, R}) = 3$. 
\end{pf}

Let $G$ be a graph and  $D$  a subset of $V(G)$ such that every vertex not in $D$ is adjacent to at least one member of $D$ then $D$ is called the dominating set for $G$. It is obvious that $V(G)$ is a dominating set for $G$. Again, it is easy to see that for any non-commutative subring $S$ of $R$, the set $S \setminus Z(S)$ is a dominating set for $\Gamma_{S, R}$. Let $A$ and $B$ be two subsets of $R$. We define $A + B := \{a + b : a \in A, b \in B\}$. Then it can be seen that $(S + C_R(S)) \setminus C_R(S)$ is a dominating set for $\Gamma_{S, R}$ if $S$ is a non-commutative subring of a finite ring $R$. 
The following proposition also gives dominating sets for $\Gamma_{S, R}$.

\begin{Prop}
Let $S$ be a subring of a ring $R$ and
$A \subseteq V(\Gamma_{S, R})$. Then $A$ is a dominating set for $\Gamma_{S, R}$ if and only if $C_R(A) \subseteq A \cup C_R(S)$. 
\end{Prop}
\begin{pf}
Suppose, $A$ is a dominating set for $\Gamma_{S, R}$ and $v \in V(\Gamma_{S, R})$ such that $v \in C_R(A)$. If $v \notin A$ then there exists an element $a \in A$ such that $va \neq av$, a contradiction.

Conversely, we suppose that $C_R(A) \subseteq A \cup C_R(S)$. Let $v \in V(\Gamma_{S, R})$ such that $v \notin A$. Suppose that $va = av$ for all $a \in A$. Then $v \in C_R(A)$ and so $v \in A \cup C_R(S)$. Thus, $v \in A$, a contradiction. Hence, $A$ is a dominating set for $\Gamma_{S, R}$.
\end{pf}

\begin{Prop}
Let $R$ be a ring with unity and $S$ a subring of $R$. If $L = \{s_1, s_2, \dots, s_n\}$ is a generating set for $S$ and $L \cap C_R(S) = \{s_{m + 1}, \dots, s_n\}$ then $K = \{s_1, s_2, \dots, s_m\} \cup \{s_1 + s_{m + 1}, s_1 + s_{m + 2}, \dots, s_1 + s_n\}$ is a dominating set for $\Gamma_{S, R}$. 
\end{Prop} 

\begin{pf}
Clearly, $K \subseteq V(\Gamma_{S, R})$. Let $v \in V(\Gamma_{S, R})$ such that $v \notin L$. If $v \in S$ then there exists an element $s =  \beta_i{s_1}^{\alpha_{1i}}{s_2}^{\alpha_{2i}}\dots  {s_d}^{\alpha_{di}}$, where $\beta_i \in \Z$, $\alpha_{ji} \in \N$ $\cup \{0\}$ and $s_j \in L$ such that $vs \neq sv$. Therefore, $vs_i \neq s_iv$ for some $1 \leq i \leq m$ and so, $v$ is adjacent to $s_i$. 

If $v \in R \setminus S$ then there exists an element $u = \gamma_i{s_1}^{\alpha_{1i}}{s_2}^{\alpha_{2i}}\dots {s_p}^{\alpha_{pi}}$, where $\gamma_i \in \Z,$ $\alpha_{li} \in \N$ $\cup \{0\}$ and $s_l \in L$ such that $vu \neq uv$. If $vs_i \neq s_iv$ for some $1 \leq i \leq m$ then $v$ is adjacent to $s_i$. Otherwise, $vs_i = s_iv$ for all $1 \leq i \leq m$. So, there exists an element $s_l$ for some $ m + 1 \leq l \leq n$ such that $vs_l \neq s_lv$. Therefore, $v$ is adjacent to $s_1 + s_l$. Hence, the proposition.
\end{pf}

An edge coloring of a graph $G$ is an assignment of ``colors" to the edges of the graph so that no two adjacent edges have the same color. The chromatic index of a graph denoted by $\chi^\prime(G)$ and is defined as the minimum number of colours needed for a colouring of $G$. Let $\vartriangle$ be the maximum vertex degree of $G$, then Vizing's theorem \cite{gt} gives $\chi^\prime(G) =  \vartriangle$ or $\vartriangle + 1$. Thus, Vizing's theorem divides the graphs into two classes according to their chromatic index. Graphs satisfying $\chi^\prime(G) = \vartriangle$ are called graphs of class $1$ and those with $\chi^\prime(G) = \vartriangle + 1$ are called graphs of class $2$. Following theorem shows that $\Gamma_{R, R}$ is of class $2$.

\begin{Thm}
Let $R$ be a ring. Then the non-commuting graph $\Gamma_{R, R}$ is of class $2$.
\end{Thm}

\begin{pf}
Clearly, $\vartriangle \leq |R| - |Z(R)| - 1$. If $\chi^\prime(\Gamma_{S, R}) = \vartriangle$ then $\chi^\prime(\Gamma_{S, R}) \leq |R| - |Z(R)| - 1$, which is not true for the ring $R = \left\lbrace  \begin{bmatrix}
    a & b\\
    0 & 0\\
  \end{bmatrix} \;:\; a, b \in Z_{2} \right\rbrace$.
Hence, $\Gamma_{R, R}$ is class of $2$.
\end{pf}

We conclude this section with the following conjecture.
\begin{Conj}
Let $S$ be a proper non-commutative subring of $R$. Then the relative non-commuting graph $\Gamma_{S, R}$ is of class $1$.
\end{Conj}

\section{Relative non-commuting graphs and $\Pr(S, R)$}
In this section, we give some connections between $\Gamma_{S, R}$ and $\Pr(S, R)$, where $S$ is a subring of a finite ring $R$. We start with the following result. 

\begin{Thm}\label{relation}
Let $S$ be a subring of a ring $R$. Then the number of edges of $\Gamma_{S, R}$ is
\begin{center}
$
|E(\Gamma_{S, R})| = |S||R|(1 - \Pr(S, R)) - \frac{|S|^2}{2}(1 - \Pr(S)). 
$
\end{center}
\end{Thm}

\begin{pf}
Let $I = \{(r_1, r_2) \in S \times R : r_1r_2 \neq r_2r_1\}$ and $J = \{(r_1, r_2) \in R \times S : r_1r_2 \neq r_2r_1\}$.   Therefore, we have 
$|I| = |S||R| - |\{(r_1, r_2) \in S \times R : r_1r_2 = r_2r_1\}| = |S||R| - |S||R|\Pr(S, R) = |J|$
and so $|I \cap J| = |\{(a, b) \in S \times S : ab \neq ba\}| = |S|^2 - |S|^2 \Pr(S).$ Thus, the result follows from the fact that  $|E(\Gamma_{S, R})| = \frac{1}{2}|I \cup J|$.
\end{pf}

The above theorem shows that lower or upper bounds for $\Pr(S)$ and $\Pr(S, R)$ will give lower or upper bounds for $|E(\Gamma_{S, R})|$ and vice-versa. More bounds for $|E(\Gamma_{S, R})|$ are obtained in the next few results.


\begin{Prop}\label{lower bound for E}
Let $S$ be a subring of a ring $R$. Then
\[
|E(\Gamma_{S, R})| \geq \frac{1}{2}|S||R| - \frac{1}{4}|S|^2 - \frac{1}{4}|Z(S)||R| - \frac{1}{4}|S||C_R(S)| + \frac{1}{4}|Z(S)||S|.
\]
\end{Prop}

\begin{pf}
Let $A = V(\Gamma_{S, R}) \cap S$ and $B = V(\Gamma_{S, R}) \cap (R \setminus S)$. Therefore, $|A| = |S| - |Z(S)|$ and $|B| = |R| - |S| - |C_R(S)| + |Z(S)|$. So, we have
\begin{align*}
2|E(\Gamma_{S, R})| &= \underset{v \in V(\Gamma_{S, R})}{\sum}\deg(v) = \underset{v \in A}{\sum}\deg(v) + \underset{v \in B}{\sum}\deg(v)\\
&= \underset{v \in A}{\sum}(|R| - |C_R(r)|) + \underset{v \in B}{\sum}(|S| - |C_S(r)|)\\ &\geq |A||R| - \frac{|A||R|}{2} - |B||S| - \frac{|S||B|}{2}.
\end{align*}
Thus, putting the values of $|A|$ and $|B|$, we get the required result.    
\end{pf}
We conclude this section with some  consequences of Theorem \ref{relation}.

\begin{Prop}
Let $S$ be a non-commutative subring of a ring $R$ and $p$ the smallest prime dividing $|R|$. Then 
\[
|E(\Gamma_{S, R})| \leq |S|(|R| - \frac{3|S|}{16} - p) - |Z(R) \cap S|(|R| - p)
\]
\end{Prop}
\begin{pf}
By \cite[Theorem 2.5]{duttabasnetnath}, we have 
\begin{equation}\label{newneweq}
\frac{|Z(R) \cap S|}{|S|} + \frac{p(|S| - |Z(R) \cap S|)}{|S||R|} \leq \Pr(S, R).
\end{equation}
Now, using \eqref{newneweq} and the fact that $\Pr(S) \leq \frac{5}{8}$  in Theorem \ref{relation} we get the required result. 
\end{pf}

\begin{Prop}
Let $S$ be a non-commutative subring of a ring $R$. Then 
\[
|E(\Gamma_{S, R})| \geq -\frac{3|S|^2}{16} + \frac{3|S||R|}{8}.
\]
\end{Prop}
\begin{pf}
Using \cite[Theorem 2.2]{duttabasnetnath}, we have that $\Pr(S, R) \leq \Pr(S) \leq \frac{5}{8}$. Therefore, $1 - \Pr(S, R) \geq 1 - \Pr(S) \geq \frac{3}{8}$. Hence, putting these results in Theorem \ref{relation}, we get the required proposition. 
\end{pf}


\begin{Prop}
Let $S$ be a non-commutative subring of a ring $R$. If $|C_R(S)| = 1$ then
\[
2|R|\Pr(S, R) - |S|\Pr(S) \neq -2\frac{|R|}{|S|} + \frac{4}{|S|} + 2|R| - |S|.
\]
\end{Prop}
\begin{pf}
Suppose there exists a finite ring $R$ with non-commutative subring $S$ such that $|C_R(S)| = 1$ and
\[
2|R|\Pr(S, R) - |S|\Pr(S) = -2\frac{|R|}{|S|} + \frac{4}{|S|} + 2|R| - |S|.
\]
Then the above equation, in view of Theorem \ref{relation}, gives
 \[
|E(\Gamma_{S, R})| = |R| - |C_R(S)| - 1 = |V(\Gamma_{S, R})| -1. 
\]
This shows that there is a finite non-commutative ring $R$ with non commutative subring $S$ such that $\Gamma_{S, R}$ is a star graph, which is not possible (by Theorem \ref{star}). Hence, the proposition follows.
\end{pf}

\section{Relative non-commuting graph and relative $\Z$-isoclinism}

In $1940$, Hall \cite{pH40} introduced the notion of isoclinism between two groups. Following Hall, Buckley et al. \cite{BMS} introduced the concept of $\Z$-isoclinism between two rings. Recently, Dutta et al. \cite{duttabasnetnath} introduced the concept of relative $\Z$-isoclinism between two pairs of rings. For a subring $S$ of $R$, $[S, R]$ is the subgroup of $(R, +)$ generated by all commutators $[s, r], s \in S, r \in R$. Let $S_1$ and $S_2$ be two subrings of the rings $R_1$ and $R_2$ respectively. Recall that a pair of rings $(S_1, R_1)$ is said to be relative $\Z$-isoclinic to a pair of rings $(S_2, R_2)$ if there exist additive group isomorphisms $\phi : \frac{R_1}{Z(R_1) \cap S_1} \rightarrow \frac{R_2}{Z(R_2) \cap S_2}$ such that $\phi(\frac{S_1}{Z(R_1) \cap S_1}) = \frac{S_2}{Z(R_2) \cap S_2}$ and $\psi : [S_1, R_1] \rightarrow [S_2, R_2]$ such that $\psi([s_1, r_1]) = [s_2, r_2]$ whenever $\phi(s_1 + (Z(R_1) \cap S_1)) = s_2 + (Z(R_2) \cap S_2)$ and $\phi(r_1 + (Z(R_1) \cap S_1)) = r_2 + (Z(R_2) \cap S_2)$ where $s_1 \in S_1, s_2 \in S_2, r_1 \in R_1, r_2 \in R_2$. Such pair of mappings $(\phi, \psi)$ is called a relative $\Z$-isoclinism from $(S_1, R_1)$ to $(S_2, R_2)$. In this section, we have the following main result.

\begin{Thm}\label{isoclinic_rncg} 
Let $S_1$ and $S_2$ be two subrings of the finite rings $R_1$ and $R_2$ respectively. Let the pairs $(S_1, R_1)$ and $(S_2, R_2)$ are relative $\Z$-isoclinic. Then $\Gamma_{S_1, R_1} \cong \Gamma_{S_2, R_2}$ if $|Z(R_1) \cap S_1| = |Z(R_2) \cap S_2|$ and $|Z(R_1)| = |Z(R_2)|$.  
\end{Thm}
\begin{pf}
Suppose $(\phi, \psi)$ is a relative $\Z$-isoclinism between 
$(S_1, R_1)$ and $(S_2, R_2)$. If $|Z(R_1) \cap S_1| = |
Z(R_2) \cap S_2|$ and $|Z(R_1)| = |Z(R_2)|$ then $|S_1| = |
S_2|, |\frac{R_1}{Z(R_1)}| = |\frac{R_2}{Z(R_2)}|, |Z(R_1) 
\setminus S_1| = |Z(R_2) \setminus S_2|$ and $|S_1 \setminus 
Z(R_1)| = |S_2 \setminus Z(R_2)|$. Now, by second isomorphism 
theorem (of groups), we have $\frac{S_1}{S_1 \cap Z(R_1)} \cong \frac{S_1 + Z(R_1)}{Z(R_1)}$. Let $\{s_1, s_2, \dots, s_m\}$ be a transversal for $\frac{S_1 + Z(R_1)}{Z(R_1)}$. So, the set $\{s_1, s_2, \dots, s_m\}$ can be extended to a transversal for $\frac{R_1}{Z(R_1)}$. Suppose, $\{s_1, s_2, \dots, s_m, r_{m + 1}, \dots, r_n\}$ is a transversal for $\frac{R_1}{Z(R_1)}$. Similarly, we can find a transversal $\{s'_1, s'_2, \dots, s'_m,\\ r'_{m + 1}, \dots, r'_n\}$ for $\frac{R_2}{Z(R_2)}$ such that $\{s'_1, s'_2, \dots, s'_m\}$ is a transversal for $\frac{S_2 + Z(R_2)}{Z(R_2)} \cong \frac{S_2}{S_2 \cap Z(R_2)}$.

Let $\phi$ be defined as $\phi(s_i + Z(R_1)) = s'_i + Z(R_2)$, $\phi(r_j + Z(R_1)) = r'_j + Z(R_2)$ for $1 \leq i \leq m, m + 1 \leq j \leq n$ and let the one-to-one correspondence $\theta : Z(R_1) \rightarrow Z(R_2)$ maps elements of $S_1$ to $S_2$. Therefore, $|C_{R_1}(S_1)| = |C_{R_2}(S_2)|$. Let us define a map $\alpha : R_1 \rightarrow R_2$ such that $\alpha(s_i + z) = s'_i + \theta(z)$, $\alpha(r_j + z) = r'_j + \theta(z)$ for $1 \leq i \leq m, m + 1 \leq j \leq n$ and $z \in Z(R_1)$. Then $\alpha$ is a bijection. This gives that $\alpha$ is also a bijection from $R_1 \setminus C_{R_1}(S_1)$ to $R_2 \setminus C_{R_2}(S_2)$. Suppose $u, v$ are adjacent in $\Gamma_{S_1, R_1}$. Then $u \in S_1$ or $v \in S_1$, say $u \in S_1$. So, $[u, v] \neq 0$, therefore $[s_i + z, r +z_1] \neq 0$, where $u = s_i + z, v = r +z_1$ for some $z, z_1 \in Z(R_1)$, $r \in \{s_1, s_2, \dots, s_m, r_{m + 1}, \dots, r_n\}$ and $1 \leq i \leq m$. Thus $[s'_i + \theta(z), r + \theta(z_1)] \neq 0$, where $\theta(z), \theta(z_1) \in Z(R_2)$ and so, $\alpha(u)$ and $\alpha(v)$ are adjacent. Hence, the theorem. 
\end{pf}

\noindent We conclude the paper with the following consequence of Theorem \ref{isoclinic_rncg}.
\begin{Cor}
Let $R$ be a ring with subrings $S$ and $T$ such that $(S, R)$ is relative $\Z$-isoclinic to $(T, R)$. Then $\Gamma_{S} \cong \Gamma_{T}$ if $|Z(R) \cap S| = |Z(R) \cap T|$.
\end{Cor}


\end{document}